\documentclass[11pt, a4paper]{amsart}
\usepackage{amsfonts}

\newtheorem{theorem}{\bf Theorem}[section]
\newtheorem{remark}{\bf Remark}[section]
\newtheorem{proposition}{Proposition}[section]
\newtheorem{lemma}{Lemma}[section]

\usepackage{enumerate,mathrsfs}
\usepackage{amssymb,amsmath,graphicx,amsthm,mathtools,hyperref}

\usepackage{booktabs}

\theoremstyle{plain}
\usepackage{siunitx}

\usepackage{algorithm}
\usepackage{algorithmic,colortbl}

\makeindex
\usepackage[intoc]{nomencl} 

\makenomenclature

\usepackage{hyperref}
\hypersetup{nesting=true,debug=true,naturalnames=true}
\usepackage{graphicx,amssymb,upref}

\usepackage{subcaption}
\usepackage{enumerate,float}
\usepackage{fullpage}
\makeatletter
\@namedef{subjclassname@2020}{\textup{2020} Mathematics Subject Classification}
\makeatother

\usepackage{lipsum}
\begin{document}

\title{ Bivariate Tempered Space-Fractional Poisson Process and Shock Models}
\author[]{Ritik Soni}
\address{\emph{Department of Mathematics and Statistics,
		Central University of Punjab, Bathinda, Punjab -151401, India.}}
\email{ritiksoni2012@gmail.com}
\author[]{Ashok Kumar pathak}
\address{\emph{Department of Mathematics and Statistics,
		Central University of Punjab, Bathinda, Punjab -151401, India.}}
\email{ashokiitb09@gmail.com}
\author[]{Antonio Di Crescenzo}
\address{\emph{Dipartimento di Matematica,
		Universit\`{a} degli Studi di Salerno, I-84084 Fisciano (SA), Italy.}}
\email{adicrescenzo@unisa.it}
\author[]{Alessandra Meoli}
\address{\emph{Dipartimento di Matematica,
		Universit\`{a} degli Studi di Salerno, I-84084 Fisciano (SA), Italy.}}
\email{ameoli@unisa.it}


			
			\begin{abstract}
				In this paper, we introduce a  bivariate tempered space-fractional Poisson process (BTSFPP)  by time-changing the bivariate Poisson process
				 with an independent tempered $\alpha$-stable subordinator. 
				 We study its distributional properties and its connection to differential equations. The  L\'{e}vy measure for the BTSFPP  is also derived. A bivariate competing risks and shock model based on the BTSFPP for predicting the failure times of the items that undergo two random shocks is  also explored.  The system is supposed to break when the sum of two types of shocks reaches a certain random threshold. Various results related to reliability such as reliability function, hazard rates, failure density, and the probability that the failure occurs due to a certain type of  shock are studied.  We show that for a general L\'{e}vy subordinator, the failure time of the  system is  exponentially distributed with mean depending on the Laplace exponent of the L\'{e}vy subordinator when the threshold has geometric distribution. Some special cases and several typical examples  are also demonstrated.
				\end{abstract}
			\maketitle
			\section{Introduction}
			The Poisson process is one of the most widely used counting processes with nice mathematical properties and applications in diverse disciplines of applied sciences, namely insurance, economics, biology, queuing theory, reliability, and statistical physics (see \cite{Byrne, Golding, Jung, Stanislavsky, Wang}). In recent years, the construction and generalization of the counting processes via subordination techniques have received a considerable amount of interest from theoretical and application view points. Orsingher and Polito \cite{Orsingher} introduced a space-fractional version of the Poisson process by subordinating the homogeneous Poisson process (HPP) with an independent $\alpha$-stable subordinator. Meerschaert et al. \cite{Meerschaert} studied the Poisson process by considering an inverse stable subordinator and established its connection with  the fractional Poisson process. Orsingher and Toaldo \cite{Orsingher1} proposed a unified approach by time-changing in HPP with an independent general L\'{e}vy subordinator. For more recent developments in this direction, one may refer to \cite{ Gajda, Kumar2, Kumar1, Maheshwari1, Maheshwari, Soni} and references therein.
			
			Apart from univariate counting processes, researchers have explored multivariate versions of the counting process in recent years for effectively analyzing  complex real-world phenomena arising in daily life.  Beghin \cite{Beghin} defined a multivariate fractional Poisson counting process by considering common random time-change of a finite-dimensional independent Poisson process. Leonenko and Merzbach \cite{Leonenko}  considered a multi-parameter fractional Poisson process using inverse subordinators and Mittag-Leffler functions and studied its main characteristics. Di Crescenzo and Longobardi \cite{Crescenzo1} discussed a bivariate Poisson process with applications in shock models. Recently, Di Crescenzo and Meoli \cite{Crescenzo} considered a bivariate space-fractional Poisson process and studied competing risks and shock models associated with it. In reliability theory and survival analysis, system failure is discussed primarily using conventional competing risks and shock models. Lehmann \cite{Lehmann} presented a class of general shock models in which failure arises as a result of a competing cause of trauma-related degradation. Cha and Giorgio \cite{Cha} developed a new class of bivariate counting processes that have
			marginal regularity property and utilizes it in a shock model. For recent development in this area, one can refer to Cha and Finkelstein \cite{Cha1}, Di Crescenzo and Meoli \cite{Crescenzo}, and Di Crescenzo and Pellerey \cite{Crescenzo2}.
			
			In this paper, we introduce  a bivariate tempered space-fractional Poisson process (BTSFPP) by time-changing the bivariate Poisson process with an independent tempered $\alpha$-stable subordinator (TSS) and study its important characteristics. We derive its L\'{e}vy  measure and the governing  
			differential equations of the probability mass function (pmf)  and probability generating function (pgf). We also propose a shock model for predicting the failure time of items subject to two external random shocks in a counting pattern governed by the BTSFPP. The system is supposed to break when two types of shocks reach their random thresholds. The results related to reliability such as  reliability function, hazard rates, failure density, and the probability that the failure occurs due to a certain type of  the shock are studied. Several typical examples based on  different random threshold distributions are also presented. Later on, for a general L\'{e}vy subordinator, we showed that the failure time of the  system is  exponentially distributed with mean depending on the Laplace exponent of the L\'{e}vy subordinator when the threshold is geometrically distributed. Graphs of survival function for different values of tempering parameters $\theta$ and stability index $\alpha$ are shown.
			
			The structure of the paper is as follows: In Section 2, we present some preliminary notations and definitions. In Section 3, we introduce the bivariate tempered space-fractional Poisson process (BTSFPP) and discuss its connection to  differential equations. A bivariate shock system governed by the BTSFPP and some results related to reliability of the failure time of the system are provided in Section 4. Also, we present a bivariate Poisson time-changed shock model when the underlying process is governed by an independent general  L\'{e}vy subordinator.  Finally, some concluding remarks are discussed in the last section.
			
			\section{Preliminaries}
		In this section, some notations and results are given which will be used in the subsequent sections.
			Let $\mathbb{N}$ denotes the set of natural numbers and $\mathbb{N}_0 = \mathbb{N} \cup \{0\}$. Let $\mathbb{R}$ and $\mathbb{C}$ denote the set of real numbers and set of complex numbers, respectively.
			\subsection{Generalized Wright Function}
			The generalized Wright function is defined by (see \cite{Kilbas})
			\begin{equation}\label{gwf11}
			_p\psi_q \left[z\; \vline \;\begin{matrix}
			\left(\alpha_i, \beta_i\right)_{1,p}\\
			(a_j,b_j)_{1,q}
			\end{matrix} \right] = \sum_{k=0}^{\infty} \frac{z^k}{k!} \frac{\prod_{i=1}^{p} \Gamma(\alpha_i + \beta_i k)}{\prod_{j=1}^{q}\Gamma(a_j + b_j k)},\;\;  z, \alpha_i, a_i \in \mathbb{C}\; \text{and}\; \beta_i, b_i \in \mathbb{R},
			\end{equation}
			under the convergence condition \begin{equation*}
			\sum_{j=1}^{q} b_j - \sum_{i=1}^{p} \beta_i >-1.
			\end{equation*}
			\subsection{L\'{e}vy Subordinator} A L\'{e}vy subordinator denoted by $\{S (t)\}_{t \geq 0}$ is a  non-decreasing L\'{e}vy process with Laplace transform (see \cite{Applebaum})
			\begin{equation*}
			\mathbb{E}\left(e^{-uS(t)}\right) = e^{-t \psi(u)},\;\; u \geq 0,
			\end{equation*}
			where $\psi(u)$ is Laplace exponent given by
			\begin{equation*}
			\psi(u) = \eta u +\int_{0}^{\infty} (1-e^{-ux}) \nu (dx), \;\; \eta \geq 0.
			\end{equation*}
			Here $\eta$ is the drift coefficient and $\nu$ is a non-negative L\'{e}vy measure on the positive half-line satisfying
			\begin{equation*}
			\int_{0}^{\infty} \min \{x,1\} \nu (dx) < \infty, \;\;\; \text{ and } \;\;\; \nu ([0,\infty)) =\infty,
			\end{equation*}
			so that $\{S (t)\}_{t \geq 0}$ has strictly increasing sample paths almost surely (a.s.).\\
			\subsubsection{\textbf{Tempered $\alpha$-Stable Subordinator}}
			For $\alpha \in (0,1)$ and $\theta >0$, the tempered $\alpha$-stable subordinator $\{S^{\alpha, \theta}(t)\}_{t\geq 0}$ is defined by the Laplace transform (see \cite{Kumar})
			\begin{equation}\label{tss11}
			\mathbb{E}[e^{-u S^{\alpha, \theta}(t) }] =  e^{\displaystyle -t\left((u+\theta)^\alpha - \theta^\alpha \right)},
			\end{equation} 
			with Laplace exponent $\psi(u) = (u+\theta)^\alpha - \theta^\alpha$.\\
			\noindent Further, the 
	L\'{e}vy measure associated with $\psi$  is (see \cite{Gupta2})
	\begin{equation}\label{lm111}
	\nu(s) = \frac{\alpha}{\Gamma(1-\alpha)}\frac{e^{-\theta s}}{s^{\alpha+1}}, \;\; s >0.
	\end{equation} 
	 Let $f_{S^{\alpha, \theta}(t)}(x,t)$ denote the probability density function (pdf) of the TSS.  By  independent and stationary increments of the L\'{e}vy subordinator, the joint density is defined as
	\begin{equation}\label{li1}
	f_{S^{\alpha, \theta}(t_1),S^{\alpha, \theta}(t_2) }(x_1, t_1; x_2,t_2)dx_1 dx_2 =  f_{S^{\alpha, \theta}(t_2-t_1) }(x_2-x_1, t_2-t_1) f_{S^{\alpha, \theta}(t_1) }(x_1, t_1)dx_1 dx_2.
	\end{equation}

	\subsection{Tempered Space-Fractional Poisson Process} Let $\{\mathcal{N}(t,\lambda)\}_{t \geq 0}$ be the homogeneous Poisson process with parameter $\lambda >0$. The tempered space-fractional Poisson process (TSFPP) denoted by $\{\mathcal{N}^{\alpha, \theta}(t,\lambda)\}_{t \geq 0}$ is defined by time-changing the homogeneous Poisson process with an independent TSS as (see \cite{Gupta})
	\begin{equation*}
	\mathcal{N}^{\alpha, \theta}(t,\lambda) := \mathcal{N}(S^{\alpha, \theta}(t), \lambda).
	\end{equation*}
	Its pmf $p^{\alpha, \theta}(k,t)$ is given by (see \cite{Gupta1})
	\begin{equation*}
	p^{\alpha, \theta}(k,t) = \frac{(-1)^k}{k!}  e^{t\theta^{\alpha}}\sum_{i=0}^{\infty} \frac{\theta^i}{\lambda^i i!} \; _1\psi_1 \left[-\lambda^\alpha t\; \vline \;\begin{matrix}
	\left(1,\alpha\right)\\
	(1-k-i, \alpha)
	\end{matrix} \right].
	\end{equation*}
	
	\subsection{Backward Shift Operators} Let $B$ is the backward shift operator defined by $B[\xi(k)] = \xi(k-1)$. For the fractional difference operator $(I-B)^\alpha$, we have  (see \cite{Orsingher})
	\begin{equation}\label{do1}
	(I-B)^\alpha = \sum_{i=0}^{\infty} \binom{\alpha}{j}(-1)^i B^i, \;\;\alpha \in (0,1),
	\end{equation}
	where $I$ is an identity operator.\\
	Furthermore, let  $\{B_i\}$, $ i\in \{1,2,\dots, m\}$ be the operators defined as 
		\begin{equation}\label{do2}
	B_i[\xi(k_1, k_2, \dots, k_m)] = \xi(k_1, k_2, \dots, k_{i}-1,\dots, k_m).
	\end{equation}
	For $m=1$ case,  $B_i$'s act same as the operator $B$.  

\section{Bivariate Tempered Space-Fractional Poisson Process}
Let $\{\mathcal{N}_i(t, \lambda_i)\}_{t\geq 0}, i=1,2$ be two independent homogeneous Poisson processes with parameter $\lambda_i, i=1,2$, respectively. Then, for  $\alpha \in (0,1)$, we define the BTSFPP $\{\mathcal{Q}^{\alpha, \theta}(t)\}_{t\geq 0}$ as
\begin{equation}\label{bp1}
\mathcal{Q}^{\alpha, \theta}(t) := \left(\mathcal{N}_1(S^{\alpha, \theta}(t), \lambda_1), \mathcal{N}_2(S^{\alpha, \theta}(t), \lambda_2)\right) :=\left( \mathcal{N}_1^{\alpha, \theta}(t,\lambda_1), \mathcal{N}_2^{\alpha, \theta}(t,\lambda_2)\right).
\end{equation}

Throughout the paper, we work with the bivariate process.
Here, we denote any arbitrary bivariate vector of constants by $\textbf{a}=(a_1, a_2)$, where  $a_1$, $a_2$ are nonnegative integers.
  Let  $\textbf{b}=(b_1, b_2)$  and $\textbf{0} = (0,0)$ is the null vector.  We write $\textbf{a} \geq  \textbf{b}$ (or $\textbf{a} \leq  \textbf{b}$)  to means that $a_i\geq b_i$ (or $a_i\leq b_i$) for $i=1,2$. Further, we denote $\textbf{k}=(k_1, k_2)$ and $\textbf{r}=(r_1, r_2)$.  

Next, we derive the pmf, pgf, and associated differential equations for the BTSFPP.
\begin{proposition}
	For $\alpha \in (0,1)$ and $\textbf{k} \geq \textbf{0}$, the pmf $q^{\alpha, \theta}(\textbf{k},t) =\mathbb{P}\{\mathcal{Q}^{\alpha, \theta}(t) =\textbf{k}\}$ is given by
	\begin{equation}\label{pmf2}
	q^{\alpha, \theta}(\textbf{k},t) = \left(-\frac{1}{\lambda_1 +\lambda_2}\right)^{k_1+k_2} \frac{\lambda_1^{k_1}\lambda_2^{k_2}}{k_1! k_2!}e^{t \theta ^\alpha}\sum_{i=0}^{\infty} \frac{\theta^i}{i! (\lambda_1 +\lambda_2)^i} \; _1\psi_1 \left[-(\lambda_1+\lambda_2)^\alpha t\; \vline \;\begin{matrix}
	\left(1,\alpha\right)\\
	(1-(k_1+k_2)-i, \alpha)
	\end{matrix} \right].
	\end{equation}
\end{proposition}
\begin{proof}
	Firstly, we have
	\begin{align}\label{pmf1}
	q^{\alpha, \theta}(\textbf{k},t) =&\; \mathbb{P}\left(\{\mathcal{Q}^{\alpha, \theta}(t) =\textbf{k}\} \cap \left\{\mathcal{N}_1^{\alpha, \theta}(t,\lambda_1) + \mathcal{N}_2^{\alpha, \theta}(t,\lambda_2) = k_1+k_2\right\}\right) \nonumber \\
	=&\; \mathbb{P}\left(\mathcal{Q}^{\alpha, \theta}(t) =\textbf{k}\; \vline \; \left\{\mathcal{N}_1^{\alpha, \theta}(t,\lambda_1) + \mathcal{N}_2^{\alpha, \theta}(t,\lambda_2) = k_1+k_2\right\}\right) \nonumber\\
	&\;\;\;\;\times \mathbb{P} \left(\mathcal{N}_1^{\alpha, \theta}(t,\lambda_1) + \mathcal{N}_2^{\alpha, \theta}(t,\lambda_2) = k_1+k_2\right).
	\end{align}

Using the conditioning argument, we get
	\begin{align*}
	\mathbb{P}\left(\mathcal{Q}^{\alpha, \theta}(t) =\textbf{k}\; \vline \; \left\{\mathcal{N}_1^{\alpha, \theta}(t,\lambda_1) + \mathcal{N}_2^{\alpha, \theta}(t,\lambda_2) = k_1+k_2\right\}\right) =  \frac{(k_1+k_2)!}{k_1! k_2!}\frac{\lambda_1^{k_1} \lambda_2^{k_2}}{(\lambda_1+\lambda_2)^{k_1+k_2}}.
	\end{align*}
	Now, we calculate\\
		$\mathbb{P} \left(\mathcal{N}_1^{\alpha, \theta}(t,\lambda_1) + \mathcal{N}_2^{\alpha, \theta}(t,\lambda_2) = k_1+k_2\right)$
	\begin{align*}
 &= \mathbb{E}\left[\mathbb{P}(N_1(r,\lambda_1) + N_2(r,\lambda_2)) = k_1 + k_2 {\Big|}_{ r = S^{\alpha, \theta}(t)} \right]\\
	&= \mathbb{E}\left[\frac{((\lambda_1+\lambda_2)r )^{k_1+k_2}}{(k_1+k_2)!}e^{-r(\lambda_1+\lambda_2)}{\Big|}_{ r = S^{\alpha, \theta}(t)} \right]\\
	&= \frac{(-1)^{k_1+k_2}}{(k_1+k_2)!}  e^{t\theta^{\alpha}}\sum_{i=0}^{\infty} \frac{\theta^i}{(\lambda_1+\lambda_2)^i i!} \; _1\psi_1 \left[-(\lambda_1+\lambda_2)^\alpha t\; \vline \;\begin{matrix}
	\left(1,\alpha\right)\\
	(1-(k_1+k_2)-i, \alpha)
	\end{matrix} \right].
	\end{align*}
	With the help of (\ref{pmf1}), we get the pmf.
\end{proof}
\begin{remark}
	When $\theta =0, $ Equation (\ref{pmf2}) reduces to the pmf of bivariate space-fractional Poisson process studied in  \cite{Crescenzo}. 
\end{remark}
\begin{theorem}
	For $\textbf{u} = (u_1,u_2) \in [0,1]^2$, the pgf  $G^{\alpha, \theta}(\textbf{u};t)$ for the BTSFPP is given by
	\begin{equation*}
	 G^{\alpha, \theta}(\textbf{u};t) = e^{\displaystyle
		-t(\left[\lambda_1(1-u_1)+ \lambda_2(1-u_2) + \theta\right]^{\alpha} - \theta^\alpha)},
	\end{equation*}
	and it satisfies the following differential equation
	\begin{equation}\label{de23}
	\frac{d}{dt} 	G^{\alpha, \theta}(\textbf{u};t) = -\left(\left[\lambda_1(1-u_1)+ \lambda_2(1-u_2) + \theta\right]^{\alpha} - \theta^\alpha\right )	G^{\alpha, \theta}(\textbf{u};t), \;\; 	G^{\alpha, \theta}(\textbf{u};0) =1.
	\end{equation}
\end{theorem}
	\begin{proof}
		For $\lambda >0$, the pgf for the TSFPP is given by (see \cite{Gupta1})
		\begin{align*}
		\mathbb{E}\left[u^{ \mathcal{N}^{\alpha, \theta}(t,\lambda)}\right] 
		&= \mathbb{E}\left[\mathbb{E}[u^{\mathcal{N}(S^{\alpha, \theta}(t), \lambda)}|S^{\alpha, \theta}(t)]\right]\\
		&= \mathbb{E}\left[e^{-\lambda (1-u)S^{\alpha, \theta}(t)}\right]\\
		&= e^{\displaystyle -t((\lambda(1-u)+\theta)^\alpha - \theta^\alpha)}.
		\end{align*}
		We define the pgf as
		\begin{equation*}
		G^{\alpha, \theta}(\textbf{u};t) =	\mathbb{E}\left[\textbf{u}^{\mathcal{Q}^{\alpha, \theta}(t)}\right] = \sum_{\textbf{k} \geq \textbf{0}}^{} u_1^{k_1} u_2^{k_2} 	q^{\alpha, \theta}(\textbf{k},t).
		\end{equation*}
		Hence, we get
		\begin{align*}\label{pgf1123}
		G^{\alpha, \theta}(u;t) =& \mathbb{E}\left[\mathbb{E}[u^{\mathcal{Q}^{\alpha, \theta}(t, \lambda)}\;|\;S^{\alpha, \theta}(t)]\right]\nonumber\\=& \mathbb{E}\left[e^{(\lambda_1 (u_1-1)+\lambda_2(u_2-1))S^{\alpha, \theta}(t)}\right] \nonumber\\=& e^{\displaystyle
			-t(\left[\lambda_1(1-u_1)+ \lambda_2(1-u_2) + \theta\right]^{\alpha} - \theta^\alpha)}.
		\end{align*}
		By calculus, we can get  (\ref{de23}) and   the condition trivially holds for $t=0$.
	\end{proof}
	\begin{theorem}
		The pmf in (\ref{pmf2}) satisfies the following differential equation
		\begin{equation*}
		\frac{d}{dt} q^{\alpha, \theta}(\textbf{k},t) = -(\lambda_1+\lambda_2)^\alpha \left( \left(I- \frac{\lambda_1 B_1+ \lambda_2 B_2 + \theta}{\lambda_1+\lambda_2}\right)^\alpha -\left( \frac{\theta}{\lambda_1+\lambda_2}\right)^\alpha\right)q^{\alpha, \theta}(\textbf{k},t), \;\;  q^{\alpha, \theta}(\textbf{0},t) = 1.
		\end{equation*} 
	\end{theorem}
	\begin{proof}
		From (\ref{de23}),  we have
		\begin{align}\label{deq11}
		\frac{d}{dt} 	G^{\alpha, \theta}(\textbf{u};t) = -(\lambda_1+\lambda_2)^\alpha \left(\left(1- \frac{\lambda_1 u_1+ \lambda_2 u_2 + \theta}{\lambda_1+\lambda_2}\right)^\alpha -\left( \frac{\theta}{\lambda_1+\lambda_2}\right)^\alpha\right)	G^{\alpha, \theta}(\textbf{u};t).
		\end{align}
		Now, we concentrate our attention to simplify the following
		\begin{align*}
		\left(1- \frac{\lambda_1 u_1+ \lambda_2 u_2 + \theta}{\lambda_1+\lambda_2}\right)^\alpha &= \left(1- \frac{ \theta}{\lambda_1+\lambda_2}- \frac{\lambda_1 u_1+ \lambda_2 u_2 }{\lambda_1+\lambda_2}\right)^\alpha \\
		&= \sum_{j \geq 0}^{} \binom{\alpha}{j} \left(1- \frac{ \theta}{\lambda_1+\lambda_2}\right)^{\alpha-j}(-1)^j \left(\frac{\lambda_1 u_1+ \lambda_2 u_2 }{\lambda_1+\lambda_2}\right)^j\\
		&=  \sum_{j \geq 0}^{} \binom{\alpha}{j} \left(1- \frac{ \theta}{\lambda_1+\lambda_2}\right)^{\alpha-j}\frac{(-1)^j}{(\lambda_1+\lambda_2)^j} \sum_{\substack{\textbf{r} \geq \textbf{0}\\ r_1+r_2 =j}}^{} \frac{j!}{r_1! r_2!} \lambda_1^{r_1}\lambda_2^{r_2} u_1^{r_1}u_2^{r_2}.
		\end{align*}
		Therefore, from (\ref{deq11}) we get\\
		$\displaystyle \frac{d}{dt} 	G^{\alpha, \theta}(\textbf{u};t)$
		\begin{align*} 
		&= -(\lambda_1+\lambda_2)^\alpha \left(\left(1- \frac{\lambda_1 u_1+ \lambda_2 u_2 + \theta}{\lambda_1+\lambda_2}\right)^\alpha \sum_{\textbf{k} \geq \textbf{0}}^{} u_1^{k_1} u_2^{k_2} 	q^{\alpha, \theta}(\textbf{k},t) -\left( \frac{\theta}{\lambda_1+\lambda_2}\right)^\alpha \sum_{\textbf{k} \geq \textbf{0}}^{} u_1^{k_1} u_2^{k_2} 	q^{\alpha, \theta}(\textbf{k},t)\right)\\
		&= -(\lambda_1+\lambda_2)^\alpha  \sum_{j \geq 0}^{} \binom{\alpha}{j} \left(1- \frac{ \theta}{\lambda_1+\lambda_2}\right)^{\alpha-j}\frac{(-1)^j}{(\lambda_1+\lambda_2)^j} \sum_{\substack{\textbf{r} \geq \textbf{0}\\ r_1+r_2 =j}}^{} \frac{j!}{r_1! r_2!} \lambda_1^{r_1}\lambda_2^{r_2}   \\ &\times \sum_{\textbf{k} \geq \textbf{0}}^{} u_1^{k_1+r_1} u_2^{k_2+r_2} 	q^{\alpha, \theta}(k,t) +(\lambda_1+\lambda_2)^\alpha \left( \frac{\theta}{\lambda_1+\lambda_2}\right)^\alpha \sum_{\textbf{k} \geq \textbf{0}}^{} u_1^{k_1} u_2^{k_2} 	q^{\alpha, \theta}(\textbf{k},t) \\
		&= -(\lambda_1+\lambda_2)^\alpha  \sum_{j \geq 0}^{} \binom{\alpha}{j} \left(1- \frac{ \theta}{\lambda_1+\lambda_2}\right)^{\alpha-j}\frac{(-1)^j}{(\lambda_1+\lambda_2)^j} \sum_{\substack{\textbf{r} \geq \textbf{0}\\ r_1+r_2 =j}}^{} \frac{j!}{r_1! r_2!} \lambda_1^{r_1}\lambda_2^{r_2}   \\ &\times \sum_{\textbf{k} \geq \textbf{r}}^{} u_1^{k_1} u_2^{k_2} 	q^{\alpha, \theta}(\textbf{k}-\textbf{r},t) +(\lambda_1+\lambda_2)^\alpha \left( \frac{\theta}{\lambda_1+\lambda_2}\right)^\alpha \sum_{\textbf{k} \geq \textbf{0}}^{} u_1^{k_1} u_2^{k_2} 	q^{\alpha, \theta}(\textbf{k},t)\\
		&= -(\lambda_1+\lambda_2)^\alpha   \sum_{\textbf{k} \geq \textbf{r}}^{} u_1^{k_1} u_2^{k_2} \sum_{j \geq 0}^{} \binom{\alpha}{j} \left(1- \frac{ \theta}{\lambda_1+\lambda_2}\right)^{\alpha-j}\frac{(-1)^j}{(\lambda_1+\lambda_2)^j} \\ &\times \sum_{\substack{\textbf{r} \geq \textbf{0}\\ r_1+r_2 =j}}^{} \frac{j!}{r_1! r_2!} \lambda_1^{r_1}\lambda_2^{r_2} 	q^{\alpha, \theta}(\textbf{k}-\textbf{r},t)   +  (\lambda_1+\lambda_2)^\alpha \left( \frac{\theta}{\lambda_1+\lambda_2}\right)^\alpha \sum_{\textbf{k} \geq \textbf{0}}^{} u_1^{k_1} u_2^{k_2} 	q^{\alpha, \theta}(\textbf{k},t).
		\end{align*}
		Since, 
		\begin{equation}\label{bso1}
		\sum_{\substack{\textbf{r} \geq \textbf{0}\\ r_1+r_2 =j}}^{} \frac{j!}{r_1! r_2!} \lambda_1^{r_1}\lambda_2^{r_2} 	q^{\alpha, \theta}(\textbf{k}-\textbf{r},t) = (\lambda_1 B_1+ \lambda_2 B_2)^j q^{\alpha, \theta}(\textbf{k},t),
		\end{equation}
	with the help of (\ref{bso1}), we can obtain the desired differential equation. 
	\end{proof}	
Next, we derive the L\'{e}vy measure for the BTSFPP.

\begin{theorem}
	The discrete L\'{e}vy measure $\mathcal{V}_{\alpha, \theta}$ for the BTSFPP is given by
	\begin{equation}\label{lm11}
	\mathcal{V}_{\alpha, \theta} (\cdot) = \sum_{k_1,k_2 >0}^{} \frac{\lambda_1^{k_1} \lambda_2^{k_2}}{k_1! k_2!}\frac{\alpha \Gamma (k_1+k_2-\alpha)}{\Gamma(1-\alpha)} \delta_{\{\textbf{k}\}}(\cdot)  (\theta +k_1 +k_2)^{\alpha -k_1-k_2},
	\end{equation}
	where $\delta_{\{\textbf{k}\}}(\cdot)$  is Dirac measure concentrated at $\textbf{k}$.
\end{theorem}
\begin{proof}
	The pmf for the bivariate Poisson process $\mathcal{N}(t) = (\mathcal{N}_1(t, \lambda_1), \mathcal{N}_2(t, \lambda_2))$ is (see \cite{Beghin})
	\begin{equation*}
	\mathbb{P}\{\mathcal{N}_1(t, \lambda_1)= k_1, \mathcal{N}_2(t, \lambda_2)= k_2 \} = \frac{\lambda_1^{k_1} \lambda_2^{k_2}}{k_1! k_2!}t^{k_1+k_2} e^{-(\lambda_{1}+\lambda_{2})t}.
	\end{equation*}
	Using (\ref{lm111}) and applying the formula from \cite[p. 197]{Ken} to calculate the L\'{e}vy measure, we get
	\begin{align*}
	\mathcal{V}_{\alpha, \theta} (\cdot)
	&= \int_{0}^{\infty} \sum_{k_1,k_2 >0}^{} \mathbb{P}\{\mathcal{N}_1(t, \lambda_1)= k_1, \mathcal{N}_2(t, \lambda_2)= k_2 \}\;  \delta_{\{\textbf{k}\}}(\cdot)\;\nu(s)\;ds\\
	&=  \sum_{k_1,k_2 >0}^{}  \frac{\lambda_1^{k_1} \lambda_2^{k_2}}{k_1! k_2!}  \delta_{\{\textbf{k}\}}(\cdot) \frac{\alpha}{\Gamma(1-\alpha)} \int_{0}^{\infty}  e^{\displaystyle-s(\theta+k_1+k_2)}s^{k_1+k_2-\alpha-1}\;ds.
	\end{align*}
	Using integral formula 3.351.3 of \cite{Gradshteyn}, we simplify as
	\begin{equation*}
	\mathcal{V}_{\alpha, \theta} (\cdot) = \sum_{k_1,k_2 >0}^{} \frac{\lambda_1^{k_1} \lambda_2^{k_2}}{k_1! k_2!}\frac{\alpha \;  (k_1+k_2-\alpha-1)!}{\Gamma(1-\alpha)} \delta_{\{\textbf{k}\}}(\cdot)  (\theta +k_1 +k_2)^{\alpha -k_1-k_2}.
	\end{equation*}
	Hence, the theorem is proved.
\end{proof}	
With the aim of calculation of hazard rates, we establish the following lemma.
\begin{lemma}\label{lm1}
	For $h \in \mathbb{N}_{0}$, we have
	\begin{align*}
	\frac{d^h}{du^h}\left[e^{-t\left((u+\theta)^\alpha-\theta^\alpha\right)}\right] =& \sum_{k=0}^{h}\frac{1}{k!} e^{-t\left((u+\theta)^\alpha-\theta^\alpha\right)}  \sum_{j=0}^{k} \binom{k}{j}t^{k}(-1)^j  \left(\left((u+\theta)^\alpha-\theta^\alpha\right)\right)^{k-j}\\
	&\times \sum_{i=0}^{j}\binom{j}{i} (\alpha i)_h
	(u+\theta)^{\alpha i - h}(-\theta^\alpha)^{j-i}.
	\end{align*}
\end{lemma}
\begin{proof}
	Let $V(u) = -t\left((u+\theta)^\alpha-\theta^\alpha\right)$ and $W(V(u)) = e^{-t\left((u+\theta)^\alpha-\theta^\alpha\right)}$. Then, applying the Hoppe's formula (see \cite{Johnson}) to the function $W(V(u))$, we get
	\begin{equation}\label{fdb1}
	\frac{d^h}{du^h} W(V(u)) = \sum_{k=0}^{h}\frac{1}{k!} e^{-t\left((u+\theta)^\alpha-\theta^\alpha\right)}T_{h,k}(V(u)),
	\end{equation}
	where $T_{h,k}(V(u))$ is computed as
	\begin{align*}
	T_{h,k}(V(u)) &= \sum_{j=0}^{k} \binom{k}{j}\left(-V(u)\right)^{k-j}\frac{d^h}{du^h}\left(V(u)\right)^j\\
	&= \sum_{j=0}^{k} \binom{k}{j}\left(t\left((u+\theta)^\alpha-\theta^\alpha\right)\right)^{k-j}\frac{d^h}{du^h}\left(-t\left((u+\theta)^\alpha-\theta^\alpha\right)\right)^j\\
	&= \sum_{j=0}^{k} \binom{k}{j}\left(t\left((u+\theta)^\alpha-\theta^\alpha\right)\right)^{k-j}(-t)^j\sum_{i=0}^{j}\binom{j}{i}(-\theta^\alpha)^{j-i} \frac{d^h}{du^h}  (u+\theta)^{\alpha i}\\
	&= \sum_{j=0}^{k} \binom{k}{j}t^{k}(-1)^j  \left(\left((u+\theta)^\alpha-\theta^\alpha\right)\right)^{k-j} \sum_{i=0}^{j}\binom{j}{i} (\alpha i)_h
	(u+\theta)^{\alpha i - h}(-\theta^\alpha)^{j-i}.
	\end{align*}
	Hence, through (\ref{fdb1}), we proved the lemma.
\end{proof}
\section{Bivariate Shock Models}

We design a shock model that is subjected to two shocks of types 1 and 2. Let $T$ be a non-negative absolutely continuous random variable which represents the failure time  of the system  subject to two possible causes of failure. Set $\zeta=n,$ which represents the failure of the system occurring due to shock of type $n$ for $ n=1,2$. We define the total number of shocks $\{\mathcal{Z}(t)\}_{t \geq 0}$ during the time interval $[0,t]$ as
\begin{equation*}
\mathcal{Z}(t)  = \mathcal{N}_1^{\alpha, \theta}(t, \lambda_1) + \mathcal{N}_2^{\alpha, \theta}(t, \lambda_1), 
\end{equation*}
where  $\mathcal{N}_1^{\alpha, \theta}(t, \lambda_1)$ and $\mathcal{N}_2^{\alpha, \theta}(t, \lambda_2)$ are processes counting the number of shocks of type $n$ for $n=1,2$, respectively during the time interval $[0,t]$.

We introduce a random threshold $L$ which takes values in the set of natural numbers. Hence, at the first time when $\mathcal{Q}(t) =L$, the failure occurs. The probability distribution and the reliability function of $L$ are respectively defined by
\begin{equation}\label{ftpmf1}
q_k = \mathbb{P}(L=k),\;\;\; k \in \mathbb{N},
\end{equation}
and 
\begin{equation*}
\overline{q}_k = \mathbb{P}(L > k), \;\;\; k \in \mathbb{N}_0.
\end{equation*}
 Let $g_{T}(t)$ be the pdf of $T$ defined as
\begin{equation*}
T = \inf \{t \geq 0: \mathcal{Z}(t) = L\}. 
\end{equation*}
Then, we have \\
\begin{equation*}
g_T(t) = g_1(t)+g_2(t),  \;\; t \geq 0,
\end{equation*}
where the sub-densities $g_n(t)$ are defined by
\begin{equation*}
g_n(t) = \frac{d}{dt}\mathbb{P}\{T \leq t, \zeta =n\},  \;n=1,2.
\end{equation*}
Also, the probability that the failure occurs due to shock of type $n$ is given by
\begin{equation}\label{zeta1}
\mathbb{P}(\zeta =n) = \int_{0}^{\infty} g_n(t)dt, \; n=1,2.
\end{equation}
Furthermore, in terms of the joint pmf, the hazard rates are given by
\begin{align}\label{abc}
h_1(k_1, k_2; t) &= \lim_{\tau \rightarrow 0^+} \frac{\mathbb{P}\left\{\mathcal{Q}^{\alpha, \theta}(t+\tau) = (k_1+1, k_2) \;\vline \;\mathcal{Q}^{\alpha, \theta}(t) = (k_1, k_2)\right\}  }{\tau},\\
h_2(k_1,k_2; t) &= \lim_{\tau \rightarrow 0^+} \frac{\mathbb{P}\left\{\mathcal{Q}^{\alpha, \theta}(t+\tau) = (k_1, k_2+1) \;\vline \;\mathcal{Q}^{\alpha, \theta}(t) = (k_1, k_2) \right\}  }{\tau}, \nonumber
\end{align}
with $(k_1, k_2) \in \mathbb{N}_0^2.$ Hence, conditioning on $L$ and with the help of (\ref{ftpmf1}), the failure densities takes the following form
\begin{equation}\label{fd1}
g_n(t) = \sum_{k=1}^{\infty} q_k \sum_{k_1 + k_2 = k-1}^{}\mathbb{P}\left\{\mathcal{Q}^{\alpha, \theta}(t) = (k_1, k_2)\right\} h_n(k_1, k_2; t),\;\; n=1,2.
\end{equation}
The reliability function of $T$ denoted by $\overline{R}_T(t) = \mathbb{P}\{T >t\}$ is given by
\begin{equation}\label{rf112}
\overline{R}_T(t) =  \sum_{k=0}^{\infty} \overline{q}_k\sum_{k_1 + k_2 = k}^{}\mathbb{P}\left\{\mathcal{Q}^{\alpha, \theta}(t) = (k_1, k_2)\right\},\;\;\;\text{with } \overline{q}_0 =1.
\end{equation}
\begin{proposition} Under the assumptions of the model in (\ref{bp1}), 
	the hazard rates $h_n(k_1, k_2; t), t \geq 0, \text{ for }  \; n=1,2$ are given by
	\begin{align}\label{hrf1}
	h_n(k_1, k_2; t) &= \alpha  \lambda_n  (\Lambda+\theta)^{\alpha-1}e^{-t (\Lambda+\theta)^\alpha}\left(\sum_{l=0}^{\infty} \frac{\theta^l}{\Lambda^l l!} \; _1\psi_1 \left[-\Lambda^\alpha t\; \vline \;\begin{matrix}
	\left(1,\alpha\right)\\
	(1-h-l, \alpha)
	\end{matrix} \right] \right)^{-1} \nonumber \\
	& \times \sum_{k=0}^{h}\frac{1}{k!} \sum_{j=0}^{k} \binom{k}{j}(-1)^j t^k \left(\left((\Lambda+\theta)^\alpha-\theta^\alpha\right)\right)^{k-j} \sum_{i=0}^{j}\binom{j}{i} (\alpha i)_h
	(\Lambda+\theta)^{\alpha i}(-\theta^\alpha)^{j-i},
	\end{align}
	where $\Lambda=\lambda_1 +\lambda_2$ and $h = k_1+k_2$.
\end{proposition}
\begin{proof}
	We fix $n=1.$ With the help of (\ref{li1}) and considering the BTSFPP as bivariate HPP with tempered $\alpha$-stable stopping time, we get
	\begin{align*}\label{112}
	\mathbb{P}&\left\{\mathcal{Q}^{\alpha, \theta}(\tau) = (k_1+1, k_2), \mathcal{Q}^{\alpha, \theta}(t) = (k_1, k_2) \right\}\\
	&= \int_{0}^{\infty} \int_{0}^{y} \mathbb{P}\left\{\mathcal{Q}^{\alpha, \theta}(y) = (k_1+1, k_2), \mathcal{Q}^{\alpha, \theta}(x) = (k_1, k_2)\right\} \\
	&\times f_{S^{\alpha, \theta}(\tau-t) }(y-x, \tau-t) f_{S^{\alpha, \theta}(t) }(x, t)\;dy dx \\
	& = \int_{0}^{\infty} \int_{0}^{y} \mathbb{P}\left\{\mathcal{N}_1(y-x, \lambda_1) =1, \mathcal{N}_2(y-x, \lambda_2)= 0\right\} \mathbb{P}\left\{\mathcal{N}_1(x, \lambda_1) = k_1, \mathcal{N}_2(x, \lambda_2)= k_2 \right\} \\
	&\times f_{S^{\alpha, \theta}(\tau-t) }(y-x, \tau-t) f_{S^{\alpha, \theta}(t) }(x, t)\;dy dx \\
	& = \int_{0}^{\infty} \int_{0}^{y} \frac{\lambda_1^{k_1+1 }\lambda_2^{k_2}}{k_1! k_2!} e^{-(\lambda_1+\lambda_2)y} x^{k_1+k_2}(y-x)f_{S^{\alpha, \theta}(\tau-t) }(y-x, \tau-t) f_{S^{\alpha, \theta}(t) }(x, t)\;dy dx.
	\end{align*}
By the change of order of integration, we get
	\begin{align*}
	\mathbb{P}&\left\{\mathcal{Q}^{\alpha, \theta}(\tau) = (k_1+1, k_2), \mathcal{Q}^{\alpha, \theta}(t) = (k_1, k_2) \right\}\\
	&=  \frac{\lambda_1^{k_1+1 }\lambda_2^{k_2}}{k_1! k_2!} \int_{0}^{\infty}\int_{x}^{\infty} x^{h} f_{S^{\alpha, \theta}(t) }(x, t) f_{S^{\alpha, \theta}(\tau-t) }(y-x, \tau-t) e^{-(\lambda_1+\lambda_2)y}(y-x)\;dydx\\
	&= \frac{\lambda_1^{k_1+1 }\lambda_2^{k_2}}{k_1! k_2!} \int_{0}^{\infty}e^{-(\lambda_1+\lambda_2)x}x^{h} f_{S^{\alpha, \theta}(t) }(x, t)\; dx\int_{0}^{\infty} ye^{-(\lambda_1+\lambda_2)y}f_{S^{\alpha, \theta}(\tau-t) }(y, \tau-t)\; dy\\
	&=  \frac{\lambda_1^{k_1+1 }\lambda_2^{k_2}}{k_1! k_2!} (-1)^{h} 	\frac{d^h}{dx^h} \mathbb{E}\left[e^{-xS^{\alpha, \theta}(t)}\right] {\Big|}_{x =\Lambda }\;\times \;	\frac{d}{dy}\left[e^{-yS^{\alpha, \theta}(\tau-t)}\right] {\Big|}_{y= \Lambda}.
	\end{align*}
	Hence, using the definition of conditional density in (\ref{abc}), the required form is obtained with the help of (\ref{pmf2}) and Lemma (\ref{lm1}).
	Also, for $n=2$ case, the proof will follow on the same line.
\end{proof}
In the next propositions, we derive the  failure densities and the reliability function of the system and obtain the distribution (\ref{zeta1}) of failure due to $n$th type of shock.
\begin{proposition}\label{prop123}
	Under the assumptions of the model in (\ref{bp1}), for $n=1,2$ and $t \geq 0$, we have the failure density of the form
	\begin{align*}
	g_n(t) =  \alpha  \lambda_n  (\Lambda+\theta)^{\alpha-1}e^{-t\left( (\Lambda+\theta)^\alpha -\theta^\alpha \right)} &\sum_{k=1}^{\infty} q_k \frac{(-1)^{k-1}}{(k-1)!}  \sum_{l=0}^{k-1}\frac{t^{l}}{l!}  \sum_{j=0}^{l} \binom{l}{j}(-1)^j  \left(\left((\Lambda+\theta)^\alpha-\theta^\alpha\right)\right)^{l-j}
	\\ 
	& \times \sum_{i=0}^{j}\binom{j}{i} (\alpha i)_{k-1}
	(\Lambda+\theta)^{\alpha i}(-\theta^\alpha)^{j-i}.
	\end{align*}
\end{proposition}
\begin{proof}
	On substituting pmf (\ref{pmf2}) and (\ref{hrf1}) to (\ref{fd1}), we get
	\begin{align*}
	g_n(t) &=  \sum_{k=1}^{\infty} q_k \sum_{k_1 + k_2 = k-1}^{} q^{\alpha, \theta}(k,t)   \;\alpha  \lambda_n  (\Lambda+\theta)^{\alpha-1}e^{-t (\Lambda+\theta)^\alpha}\left(\sum_{l=0}^{\infty} \frac{\theta^l}{\Lambda^l l!} \; _1\psi_1 \left[-\Lambda^\alpha t\; \vline \;\begin{matrix}
	\left(1,\alpha\right)\\
	(1-h-l, \alpha)
	\end{matrix} \right] \right)^{-1} \nonumber \\
	& \times \sum_{l=0}^{k_1+k_2}\frac{1}{l!} \sum_{j=0}^{l} \binom{l}{j}(-1)^j t^l \left(\left((\Lambda+\theta)^\alpha-\theta^\alpha\right)\right)^{l-j} \sum_{i=0}^{j}\binom{j}{i} (\alpha i)_{k_1+k_2}
	(\Lambda+\theta)^{\alpha i}(-\theta^\alpha)^{j-i}\\
	&= \alpha  \lambda_n  (\Lambda+\theta)^{\alpha-1}e^{-t\left( (\Lambda+\theta)^\alpha -\theta^\alpha \right)} \sum_{k=1}^{\infty} q_k \frac{(-1)^{k-1}}{(k-1)! \Lambda^{k-1}} \sum_{k_1 =0}^{k-1} \frac{(k-1)! \lambda_{1}^{k_1} \lambda_{2 }^{k-1-k_1}}{k_1! (k-1-k_1)!} \\
	& \times \sum_{l=0}^{k-1}\frac{t^l}{l!} \sum_{j=0}^{l} \binom{l}{j}(-1)^j  \left(\left((\Lambda+\theta)^\alpha-\theta^\alpha\right)\right)^{l-j} \sum_{i=0}^{j}\binom{j}{i} (\alpha i)_{k-1}
	(\Lambda+\theta)^{\alpha i}(-\theta^\alpha)^{j-i}.
	\end{align*}
Using the binomial theorem, the failure density is obtained.
\end{proof}

\begin{proposition}\label{rf1}
	Under the assumptions of the model in (\ref{bp1}), the reliability function of $T$ is given by
	\begin{equation}\label{rf11}
	\overline{R}_T(t) =  \sum_{k=0}^{\infty} \overline{q}_k \frac{(-1)^k}{k!} \;e^{t\theta^{\alpha}} \sum_{i=0}^{\infty} \frac{\theta^i}{\Lambda^i i!} \; _1\psi_1 \left[-\Lambda^\alpha t\; \vline \;\begin{matrix}
	\left(1,\alpha\right)\\
	(1-k-i, \alpha)
	\end{matrix} \right], \;\;t \geq 0.
	\end{equation}
\end{proposition}

\begin{proof}
	The reliability function (\ref{rf11}) can be obtained  by substituting (\ref{pmf2}) to (\ref{rf112}) and simplifying it using the binomial theorem as carried out in previous  proof.
\end{proof}

\begin{proposition} Under the assumptions of the model in (\ref{bp1}),
	for $n=1,2$, we also have
	\begin{align*}
	\mathbb{P}(\zeta =n) = \alpha  \lambda_n \frac{(\Lambda+ \theta)^{\alpha-1}}{(\Lambda+\theta)^\alpha-\theta^\alpha} &\sum_{k=1}^{\infty} q_k \frac{(-1)^{k-1}}{(k-1)!}  \sum_{l=0}^{k-1}  \sum_{j=0}^{l} \binom{l}{j}(-1)^j  \left((\Lambda+\theta)^\alpha-\theta^\alpha\right)^{-j}\\& \times \sum_{i=0}^{j}\binom{j}{i} (\alpha i)_{k-1}
	(\Lambda+\theta)^{\alpha i}(-\theta^\alpha)^{j-i}.
	\end{align*}
\end{proposition}

\begin{proof}
	With the help of Proposition (\ref{prop123}), the probability (\ref{zeta1}) gives
	\begin{align*}
	\mathbb{P}(\zeta =n) = \alpha  \lambda_n  (\Lambda+\theta)^{\alpha-1} &\sum_{k=1}^{\infty} q_k \frac{(-1)^{k-1}}{(k-1)!}  \sum_{l=0}^{k-1}\frac{1}{l!}  \sum_{j=0}^{l} \binom{l}{j}(-1)^j  \left(\left((\Lambda+\theta)^\alpha-\theta^\alpha\right)\right)^{l-j}
	\\ 
	& \times \sum_{i=0}^{j}\binom{j}{i} (\alpha i)_{k-1}
	(\Lambda+\theta)^{\alpha i}(-\theta^\alpha)^{j-i} \int_{0}^{\infty}e^{-t\left( (\Lambda+\theta)^\alpha -\theta^\alpha \right)}t^{l}dt.
	\end{align*}
	Using integral formula 3.351.3 of \cite{Gradshteyn}, we get the proposition.
\end{proof}

\subsection{Generalized Shock Models}
Let $S:= \{S\left(t\right)\}_{t \geq 0}$ be a L\'{e}vy subordinator.
%
%
In the next theorem, we evaluate the reliability function of $T$ when the threshold $L$ has  geometric distribution with parameter $p\in\left(0;1\right]$, {\it i.e}.,
\begin{equation}\label{SurGeo}
\overline{q}_k=\left(1-p\right)^k, \qquad k=0,1,2\dots,
\end{equation}

\noindent and when the shocks arrive according to a process $N:=\{N\left(t\right)\}_{t \geq 0}$, where 

\begin{equation}\label{gm11}
N\left(t\right):=\left(\mathcal{N}_1\left(S\left(t\right), \lambda_{1}\right),\mathcal{N}_2\left(S\left(t\right), \lambda_{2}\right)\right),
\end{equation}

\noindent where the components of $N$ are two time-changed independent homogeneous Poisson processes with intensities $\lambda_1>0$ and $\lambda_2>0$, respectively. The time-change is represented by an independent generic subordinator $S$. 


\begin{theorem}\label{Thm1}
	For $(x_1, x_2) \in \mathbb{N}_{0}^2$ and under the assumptions of the model in (\ref{SurGeo}) and (\ref{gm11}), we have the reliability function of T as
	\begin{equation}\label{eq24}
	\overline{F}_{T}\left(t\right) = e^{-t\psi\left(\left(\lambda_1+\lambda_2\right)p\right)},
	\end{equation}
	where $\psi(\cdot)$ is the Laplace exponent of the subordinator $S$.
\end{theorem}
\begin{proof}
	Consider the reliability function of $T$ as
	\begin{align*}
	\overline{F}_{T}\left(t\right)&=\sum_{k=0}^{+\infty}\left(1-p\right)^k\sum_{x_1=0}^{k}\mathbb{P}\left(\mathcal{N}_1\left(S\left(t\right), \lambda_{1}\right)=x_1,\mathcal{N}_2\left(S\left(t\right), \lambda_{2}\right)=k-x_1\right)\\
	&=\sum_{k=0}^{+\infty}\left(1-p\right)^k\sum_{x_1=0}^{k}\frac{\lambda_{1}^{x_{1}}}{x_1!}\frac{\lambda_{2}^{k-x_{1}}}{\left(k-x_2\right)!}\int_{0}^{+\infty}e^{-\left ( \lambda _{1}+\lambda_{2} \right )s}s^{k}\mathbb{P}\left(S\left(t\right)\in\mathrm{d}s\right).
	\end{align*}
	
	\noindent We exchange the order of summation and rearrange all the terms, so to get:
	
	\begin{align}\label{SurT}
	\overline{F}_{T}\left(t\right)&=\sum_{x_1=0}^{+\infty}\frac{\lambda_{1}^{x_1}\left(1-p\right)^{x_1}}{x_{1}!}\sum_{h=0}^{+\infty}\frac{\left [ \lambda_{2} \left ( 1-p \right )\right ]^{h}}{h!}\int_{0}^{+\infty}s^{x_1+h}e^{-\left ( \lambda _{1}+\lambda_{2} \right )s}\mathbb{P}\left(S\left(t\right)\in\mathrm{d}s\right)\nonumber\\
	&=\sum_{x_1=0}^{+\infty}\frac{\lambda_{1}^{x_1}\left(1-p\right)^{x_1}}{x_{1}!}\int_{0}^{+\infty}s^{x_1}e^{-\left ( \lambda _{1}+\lambda_{2} \right )s+\lambda_2\left(1-p\right)s}\mathbb{P}\left(S\left(t\right)\in\mathrm{d}s\right)\nonumber\\
	&=\int_{0}^{+\infty}e^{-\left ( \lambda _{1}+\lambda_{2} \right )s+\lambda_2\left(1-p\right)s+\lambda_1\left(1-p\right)s}\mathbb{P}\left(S\left(t\right)\in\mathrm{d}s\right)\nonumber\\
	&=\int_{0}^{+\infty}e^{-\left ( \lambda _{1}+\lambda _{2} \right )ps}\mathbb{P}\left(S\left(t\right)\in\mathrm{d}s\right)\nonumber\\
	&=e^{-t\psi\left(\left(\lambda_1+\lambda_2\right)p\right)}.
	\end{align}
	Hence, the theorem is proved.
\end{proof}
\begin{remark}
	 In Equation (\ref{eq24}), it is observed that the random failure time $T$ is exponentially distributed with mean depending on the Laplace exponent of the subordinator. 
\end{remark}

\begin{remark}\label{rem1}
		 As a corollary of the Theorem \ref{Thm1}, it is straightforward to show that if the distribution of $L$ is a mixture of the geometric distribution (\ref{SurGeo}), then the distribution of $T$ is a mixture of the exponential distribution (\ref{SurT}). That is,
	
	\begin{equation*}
	\overline{F}_{T}\left(t\right)=\int_{0}^{1}e^{-t\psi\left(\left(\lambda_1+\lambda_2\right)p\right)}\mathrm{d}G\left(p\right),
	\end{equation*}
	
	\noindent where $G$ is a distribution on $(0,1)$.
\end{remark}
Next, we discuss examples of some special random thresholds under the assumptions of the model in (\ref{bp1}).
\subsection{Some Examples}

First, we reproduce the following identity from \cite{Gupta} as
\begin{equation}\label{tsf1}
\exp\left[-t\left(\left(\Lambda(1-u)+\theta\right)^\alpha - \theta^\alpha\right) \right] =  e^{t\theta^{\alpha}} \sum_{k=0}^{\infty} \frac{(-u)^k}{k!} \sum_{i=0}^{\infty} \frac{\theta^i}{\Lambda^i i!} \; _1\psi_1 \left[-\Lambda^\alpha t\; \vline \;\begin{matrix}
\left(1,\alpha\right)\\
(1-k-i, \alpha)
\end{matrix} \right].
\end{equation} 
Now, we derive the reliability function of $T$ for some particular cases of the random threshold $L$.\\
(I) Let $L$ follow the discrete exponential distribution with reliability function 
\begin{equation*}
\overline{q}_k = e^{-k},\; k=0,1,2,\ldots.
\end{equation*} 
From (\ref{rf1}) and with help of (\ref{tsf1}), we get
\begin{equation*}
\overline{R}_T(t) = \exp\left[-t\left(\left(\Lambda(1-e^{-1})+\theta\right)^\alpha - \theta^\alpha\right) \right].
\end{equation*}
Also, the density is given by
\begin{equation*}
g_T(t) = \frac{d}{dt} \overline{R}_T(t) = \left(\left(\Lambda(1-e^{-1})+\theta\right)^\alpha - \theta^\alpha\right)\exp\left[-t\left(\left(\Lambda(1-e^{-1})+\theta\right)^\alpha - \theta^\alpha\right) \right].
\end{equation*}
Therefore, the hazard rate function denoted by $H_T(t)$ for the random variable $T$ is given by
\begin{equation*}
H_T(t) = \frac{	g_T(t)}{	\overline{R}_T(t)} = \left(\left(\Lambda(1-e^{-1})+\theta\right)^\alpha - \theta^\alpha\right),\;\; t \geq 0.
\end{equation*} 
(II) Let $L$ follow the Yule-Simon distribution with parameter $p$ and the reliability function
\begin{equation*}
\overline{q}_k = kB(k,p+1), \; k=1,2,\ldots,
\end{equation*}
where $B(a,b) = \int_{0}^{1} t^{a-1}(1-t)^{b-1}dt$ is the beta function.
Then, the reliability function $	\overline{R}_T(t)$ takes the form
\begin{align*}
\overline{R}_T(t) &= \sum_{k=0}^{\infty} kB(k,p+1) \frac{(-1)^k}{k!} \;e^{t\theta^{\alpha}} \sum_{i=0}^{\infty} \frac{\theta^i}{\Lambda^i i!} \; _1\psi_1 \left[-\Lambda^\alpha t\; \vline \;\begin{matrix}
\left(1,\alpha\right)\\
(1-k-i, \alpha)
\end{matrix} \right]\\
&= \sum_{k=0}^{\infty} k\left(\int_{0}^{1} z^{k-1}(1-z)^{p}dz\right) \frac{(-1)^k}{k!} \;e^{t\theta^{\alpha}} \sum_{i=0}^{\infty} \frac{\theta^i}{\Lambda^i i!} \; _1\psi_1 \left[-\Lambda^\alpha t\; \vline \;\begin{matrix}
\left(1,\alpha\right)\\
(1-k-i, \alpha)
\end{matrix} \right]\\
&= e^{t\theta^{\alpha}} \int_{0}^{1}(1-z)^{p} \sum_{k=1}^{\infty} k z^{k-1}  \frac{(-1)^k}{k!}  \sum_{i=0}^{\infty} \frac{\theta^i}{\Lambda^i i!} \; _1\psi_1 \left[-\Lambda^\alpha t\; \vline \;\begin{matrix}
\left(1,\alpha\right)\\
(1-k-i, \alpha)
\end{matrix} \right]dz\\
&= \int_{0}^{1}(1-z)^{p} \left(\frac{d}{dz} \exp\left[-t\left(\left(\Lambda(1-z)+\theta\right)^\alpha - \theta^\alpha\right) \right]\right)dz.
\end{align*}
Hence, the density function is given by 
\begin{equation*}
g_T(t) =  \int_{0}^{1}(1-z)^{p} \frac{d}{dt} \left(\frac{d}{dz} \exp\left[-t\left(\left(\Lambda(1-z)+\theta\right)^\alpha - \theta^\alpha\right) \right]\right)dz.
\end{equation*}
Considering these $g_T(t)$ and $\overline{R}_T(t)$, we  get  the hazard rate function in case of Yule-Simon threshold.

Now, we discuss some special cases of the mixing distribution from Remark \ref{rem1}  under the assumptions of the model in (\ref{gm11}).

\subsection{Special Cases}

\noindent We now analyze three special cases by specifying the mixing distribution, under the assumption that $S$ is the tempered $\alpha$-stable subordinator as in (\ref{tss11}). The evaluation of the reliability functions are performed using Mathematica.\\
(I) $\mathrm{d}G\left(p\right)=\mathrm{d}p$ (uniform distribution)\\ It is 
	\begin{equation}\label{Uniform}
	\overline{F}_{T}\left(t\right) = \frac{e^{t\theta^{\alpha }}}{\alpha\left(\lambda_1+\lambda_2\right )}\times\left [ \theta E_{\frac{\alpha -1}{\alpha}}\left ( t\theta ^{\alpha } \right )-\left ( \lambda_{1}+\lambda_{2}+\theta  \right )E_{\frac{\alpha -1}{\alpha}}\left ( t\left ( \lambda_{1}+\lambda_{2}+\theta  \right )^{\alpha } \right )\right ], 
	\end{equation}
	where $E_{l}\left(z\right)=\int_{1}^{+\infty}\frac{e^{-uz}}{u^{l}}\mathrm{d}u$ is a generalized exponential integral.
	\vspace{.2cm}
\\
\noindent (II) $\mathrm{d}G\left(p\right)=\displaystyle \frac{ab\left ( 1+ap \right )^{-\left ( b+1 \right )}}{1-\left ( 1+a \right )^{-b}}\mathrm{d}p$, with $a>0$ and $b>-1\, \wedge \,b\neq 0$. (truncated Lomax). \\
	Set $a:=\frac{\lambda _{1}+\lambda _{2}}{\theta }$ and $b+1:=\alpha$. It is
	\begin{equation}\label{Lomax}
	\overline{F}_{T}\left(t\right)=\frac{e^{t\theta ^{\alpha}}\left ( \alpha -1 \right )}{\alpha \left [1-\left ( 1+\frac{\lambda _{1}+\lambda _{2}}{\theta } \right )^{1-\alpha }  \right ]}
	\times\left [ E_{2-\frac{1}{\alpha }} \left ( t\theta ^{\alpha} \right )-\left (1+ \frac{\lambda _{1}+\lambda _{2} }{\theta } \right )^{1-\alpha }E_{2-\frac{1}{\alpha }}\left (  t\theta ^{\alpha}\left ( 1+\frac{\lambda _{1} +\lambda _{2}}{\theta} \right )^{\alpha }\right )\right ].
	\end{equation}
(III) $\mathrm{d}G\left(p\right)=\frac{b}{a}\left(\frac{p-c}{a}\right)^{b-1}e^{-\left(\frac{p-c}{a}\right)^{b}}\mathrm{d}p$, where $a$ and $b$ are positive values, and $c$ is a real value. (truncated three-parameter Weibull).\\
	Set $a=\frac{1}{\lambda_1+\lambda_2}$, $b=\alpha$ and $c=-\frac{\theta}{\lambda_1+\lambda_2}$. It is
	\begin{equation}\label{Weibull}
	\overline{F}_{T}\left(t\right)=\frac{1-e^{-\left ( t+1 \right )\left [ \left ( \lambda _{1} +\lambda _{2 }+\theta \right ) ^{\alpha }-\theta ^{\alpha }\right ]}}{\left ( t+1 \right )\left [ 1-e^{-\left [ \left ( \lambda _{1} +\lambda _{2 }+\theta \right ) ^{\alpha }-\theta ^{\alpha }\right ]} \right ]}.
	\end{equation}
	
The graphs in Figures 1,2 and 3 illustrate the special cases for some particular values of parameters.
\begin{figure}[h]
	\begin{center}
		\includegraphics[scale=.52]{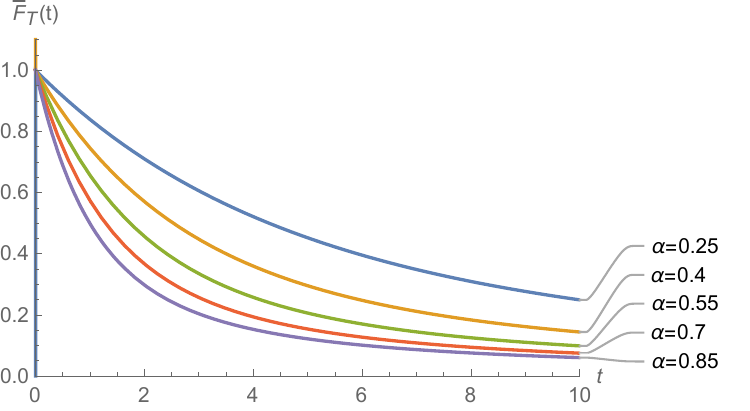}\quad
		\includegraphics[scale=.52]{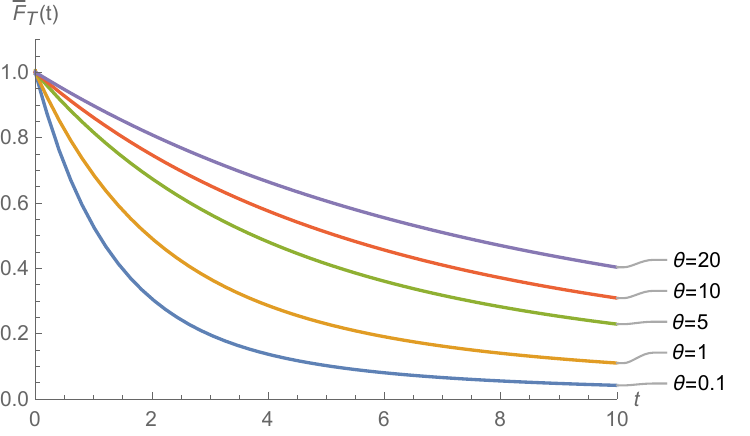}
	\end{center}
	\vspace{0.3cm}
	\caption{ Plots of the reliability function (\ref{Uniform}) with $\lambda_1=\lambda_2=1$ and $\theta=1$ on the left-hand side, $\alpha=0.5$ on the right-hand side.}
\end{figure}

\begin{figure}[h]
	\begin{center}
		\includegraphics[scale=.52]{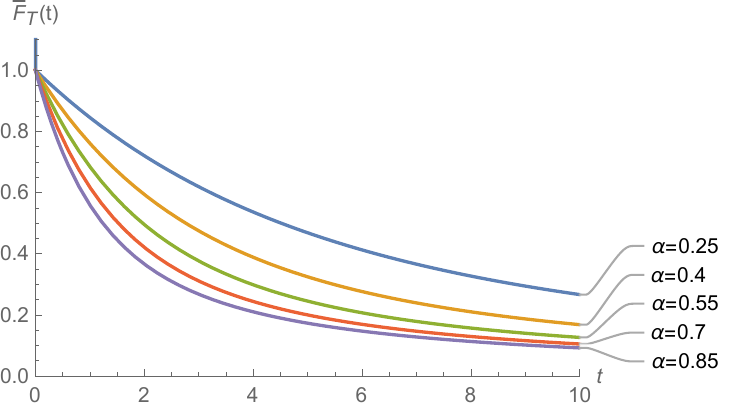}\quad
		\includegraphics[scale=.52]{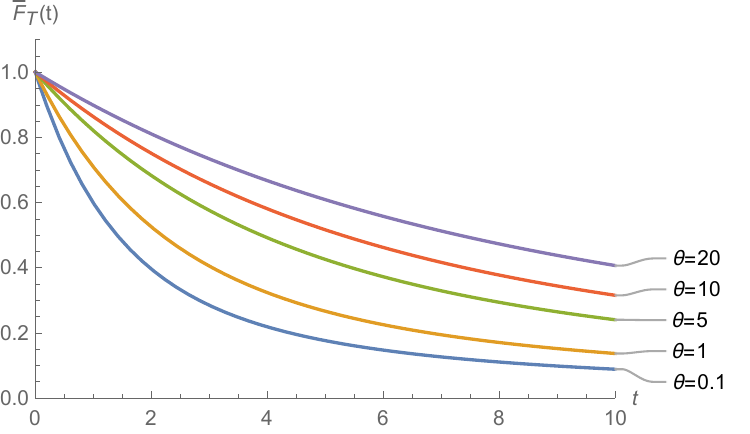}
	\end{center}
	\vspace{0.3cm}
	\caption{ Plots of the reliability function (\ref{Lomax}) with $\lambda_1=\lambda_2=1$ and $\theta=1$ on the left-hand side, $\alpha=0.5$ on the right-hand side.}
\end{figure}

\begin{figure}[h]
	\begin{center}
		\includegraphics[scale=.52]{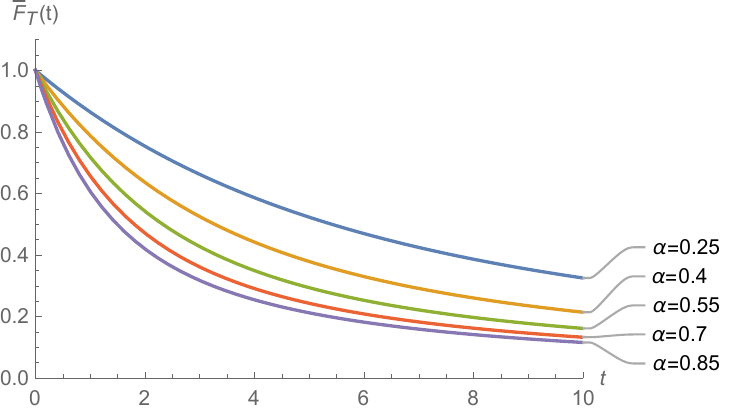}\quad
		\includegraphics[scale=.52]{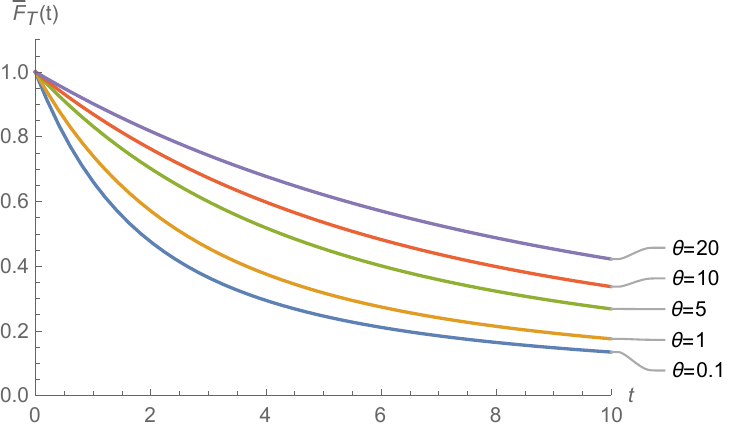}
	\end{center}
	\vspace{0.3cm}
	\caption{ Plots of the reliability function (\ref{Weibull}) with $\lambda_1=\lambda_2=1$ and $\theta=1$ on the left-hand side, $\alpha=0.5$ on the right-hand side.}
\end{figure}

\section{Concluding Remarks}
\noindent In this paper, we have proposed a bivariate tempered space-fractional Poisson process (BTSFPP) by time-changing the bivariate Poisson process with an independent tempered $\alpha$-stable subordinator. First we derived the expression for the probability mass function  and expressed it in terms of the generalized Wright function, then we obtained the governing differential equations for the pmf and the pgf. We also derived the  L\'{e}vy measure density for the BTSFPP. Influenced by the reliability applications, we  presented a bivariate competing risks and shock model based on the BTSFPP and derived various  reliability quantities to predict the life of the system. Finally, we discussed a generalized shock model and various typical examples.
	\subsection*{Acknowledgements} 
	
	
	A.D.C. and A.M. are members of the group GNCS of INdAM (Istituto Nazionale di Alta Matematica).
	
	
	\subsection*{Funding} 
	
	
	This work is partially supported by CSIR India (File No: 09/1051(11349)/2021-EMR-I), DST-SERB and MIUR--PRIN 2017, Project Stochastic Models for Complex Systems (no. 2017JFFHSH).

	\subsection*{Competing Interests} 
	
	There were no competing interests to declare which arose during the preparation or publication process of this article.
			
		\end{document}